\documentclass[conference]{IEEEtran}
\usepackage{times}

\usepackage{graphicx}
\usepackage{algorithm}
\usepackage{algorithmicx}
\usepackage{algpseudocode}
\usepackage{amsmath}
\usepackage{amssymb}
\usepackage{optidef}
\usepackage{pgfplots}
\pgfplotsset{compat=1.18} 
\usepackage{tikz}
\usepackage{tikzscale}  
\usepackage{pgf-pie}   
\usepackage{tikz-3dplot}
\usetikzlibrary{patterns}
\usepgfplotslibrary{groupplots}
\usepackage{adjustbox}
\usepackage{pifont}
\usepackage[svgnames]{xcolor}

\newcommand{\cmark}{\textcolor{ForestGreen}{\ding{51}}}
\newcommand{\xmark}{\textcolor{Red}{\ding{55}}}
\usepackage{booktabs}
\usepackage{subfig} 
\usepackage{subcaption}

\usepackage{bm}
\usepackage{comment}

\makeatletter
\renewcommand{\subsubsection}[1]{\@startsection{subsubsection}{3}{\z@}%
  {-3.25ex\@plus -1ex \@minus -.2ex}%
  {1.5ex \@plus .2ex}%
  {\normalfont\normalsize\itshape}%
  {#1}}
\makeatother

\newcommand{\R}[1]{\mathbf{R}^{#1}}
\newcommand{\jac}[2]{\frac{\partial { #1  }}{ \partial { #2 }}}

\usepackage[numbers]{natbib}
\usepackage{multicol}
\usepackage[bookmarks=true]{hyperref}

\begin{document}

\title{The Trajectory Bundle Method: \\Unifying Sequential-Convex Programming and Sampling-Based Trajectory Optimization}

\author{
Kevin Tracy$^{1}$, John Z. Zhang$^{1}$, Jon Arrizabalaga$^{1}$, Stefan Schaal$^{2}$, \\Yuval Tassa$^{3}$, Tom Erez$^{3}$,
and Zachary Manchester$^{1}$
\\
\\
$^{1}$CMU,
$^{2}$Google Intrinsic,
$^{3}$Google DeepMind\\
}



%

\maketitle

\begin{abstract}
We present a unified framework for solving trajectory optimization problems in a derivative-free manner through the use of sequential convex programming. Traditionally, nonconvex optimization problems are solved by forming and solving a sequence of convex optimization problems, where the cost and constraint functions are approximated locally through Taylor series expansions. This presents a challenge for functions where differentiation is expensive or unavailable. 
In this work, we present a derivative-free approach to form these convex approximations by computing samples of the dynamics, cost, and constraint functions and letting the solver interpolate between them. Our framework includes sample-based trajectory optimization techniques like model-predictive path integral (MPPI) control as a special case and generalizes them to enable features like multiple shooting and general equality and inequality constraints that are traditionally associated with derivative-based sequential convex programming methods. The resulting framework is simple, flexible, and capable of solving a wide variety of practical motion planning and control problems. 
\end{abstract}

\IEEEpeerreviewmaketitle

\section{Introduction}\label{sec:bundles:introduction}
Linear dynamical systems of the form $x_{+} = Ax + Bu$ underpin many of the foundational methods in modern optimal control. Ideas such as the Linear-Quadratic Regulator (LQR) and convex trajectory optimization can reason about dynamical systems of this form in a way that is globally optimal \cite{kalman1960,borrelli2017,boyd2004}. As a result, these techniques are often applied to nonlinear systems where the dynamics are locally approximated as linear around a linearization point \cite{slotine1991}. In many cases, this approximation is appropriate given the function is not being evaluated too far from where the approximation was formed.  When used appropriately, this method of linearizing nonlinear systems can be extremely effective in practice, even for highly nonlinear systems.  The two caveats here are that the nonlinear system must be both smooth and differentiable. 

For many systems, such as robotic arms, quadrotors, and wheeled vehicles, this assumption of smooth differentiability is reasonable. For rigid-body dynamics, there are specialized methods for computing derivatives of the continuous-time dynamics in a fast and efficient way \cite{featherstone1987}. However, for more complex dynamics models, there are scenarios where these derivatives are unavailable, prohibitively expensive to compute, or unreliable. If the dynamics model is learned from data, the approximation of the dynamics function may be good, while the approximation of the derivatives may be very poor. This scenario is often explored in the context of model-predictive path-integral (MPPI) control, where a learned simulator is only used to produce parallelized simulations \cite{williams2016,wagener2019}.

Another scenario in which derivatives are unavailable or unusable is in the presence of systems that make or break contact. While there has been a lot of recent interest in making contact simulation differentiable \cite{freeman2021,newbury2024,pang2023,tracy2023b,suh2022a,howell2022}, there remains a strong need for optimal control methods that do not rely on these derivatives at all.  
\begin{figure}
    \centering
    \includegraphics[width=0.95\linewidth]{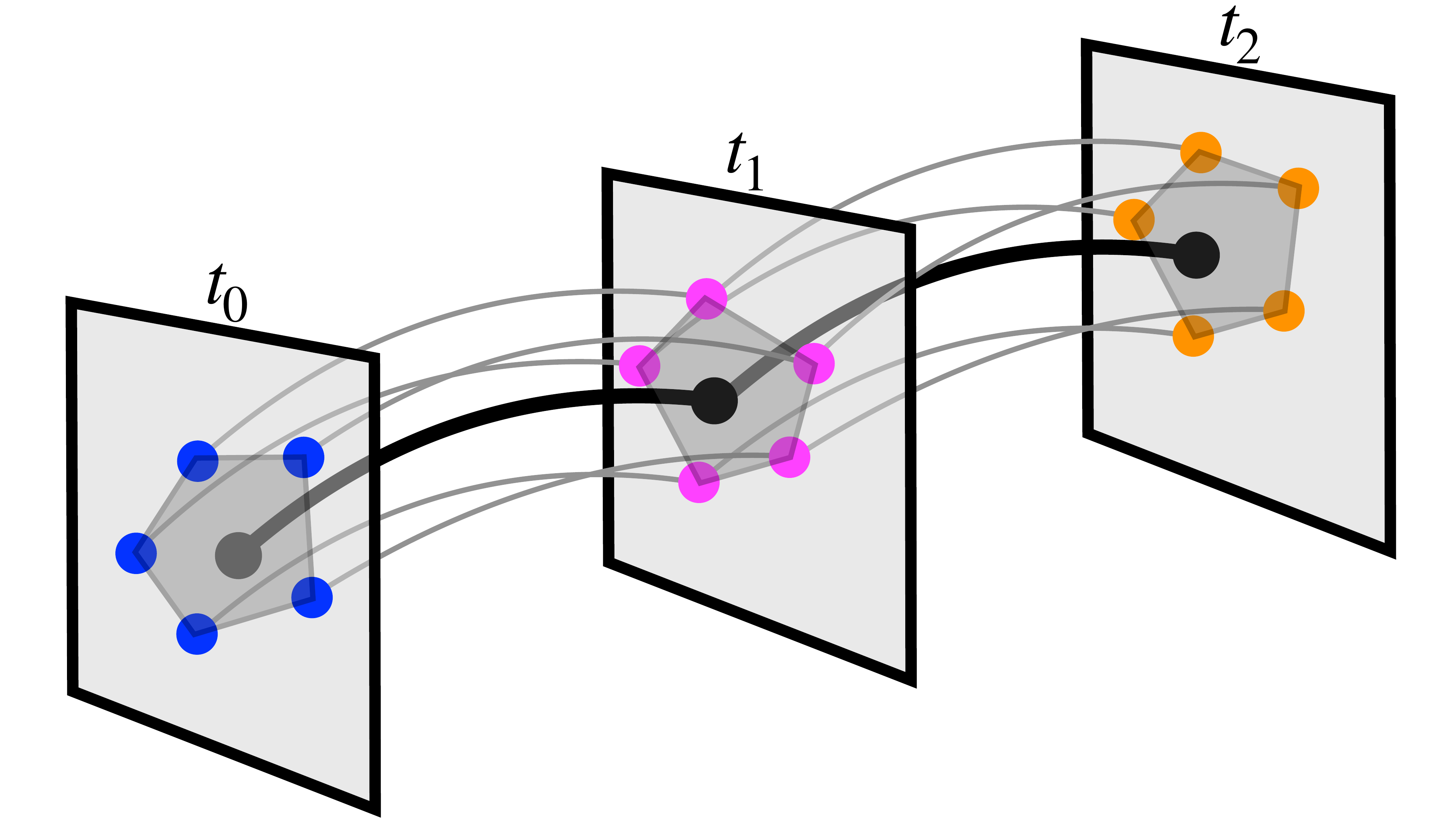}
    \caption{The Trajectory Bundle Method (TBM) is a derivative-free framework capable of solving both single-shooting sampling-based MPC problems \emph{and} general multiple-shooting trajectory optimization problems.
    TBM finds the optimal state and control sequences (bolded black trajectory) by computing samples (colored points) at each knot point ($t$s) around the current solution and using convex optimization to linearly interpolate between the samples (dark gray regions).}
    \label{fig:trajectory_bundle_method}
\end{figure}

A recent trend in robotic simulation is the introduction of simulators that can be run on accelerators for massively parallel simulation. Popular simulators like Isaac Sim \cite{makoviychuk2021,mittal2023}, Brax \cite{freeman2021}, and MuJoCo XLA (MJX) \cite{todorov2012}, are all capable of running thousands of simulations in parallel. This paper leverages the innovations in parallel simulation to motivate a new derivative-free optimal control paradigm where simulation rollouts are used to fully describe the dynamics and cost landscapes present in the problem. We introduce the trajectory bundle method for solving nonconvex trajectory optimization problems, which uses interpolated trajectories instead of derivative-based linearizations to approximate the cost, dynamics, and constraint functions in the problem. The result is a simple and robust trajectory optimization framework that can fully utilize parallelized simulation without requiring any derivatives. 

Our specific contributions in this paper are the following:
\begin{itemize}
    \item A unified framework, which we refer to as the trajectory bundle method, for solving general trajectory optimization problems using derivative-free sequential-convex programming.
    \item A method for approximating general nonlinear or non-convex cost and constraint functions through sampling and linear interpolation.
    \item A set of numerical experiments demonstrating the effectiveness of the trajectory bundle method and its equivalence to its SCP and MPPI counterparts. 
\end{itemize}


The remainder of the paper is organized as follows: we first review related literature on derivative-free optimization, sequential-convex programming, and MPPI in Section \ref{sec:related_work}. Next, we introduce relevant background on affine function approximation and its application to constrained optimization in Section \ref{sec:bundles:background}. In Section \ref{sec:trajectory_bundle_method}, we describe the trajectory bundle method in a general multiple-shooting framework and a single-shooting special case that is equivalent to MPPI. Finally, we present an array of numerical experiments in Section \ref{sec:experiments} and point avenues of future research in \ref{sec:conclusions}.

\section{Related Work}\label{sec:related_work}
In order to contextualize the trajectory bundle method, this section provides brief overviews of derivative-free optimization, trajectory optimization through sequential convex programming, and model predictive path integral control. 
\subsection{Derivative-Free Optimization}
Derivative-Free Optimization (DFO) is a well-studied technique for solving optimization problems where derivatives are unavailable. John Dennis describes the DFO problem in \cite{powell1994} as ``finding the deepest point of a muddy lake,
given a boat and a plumb line, when there is a price to be paid for each sounding.'' With modern accelerators capable of massively parallel dynamics, cost, and constraint evaluation, there is still a price to take soundings, but we can now ``buy in bulk.'' In a seminal 1965 paper, the Nelder-Mead method for derivative-free function minimization was proposed, where a simplex of sample points is used to approximate the cost landscape instead of derivatives. Methods like Mesh Adaptive Direct Search (MADS) \cite{audet2006} and NOMAD \cite{ledigabel2011,audet2021} followed with improvements to the Nelder-Mead method. 


Although there has been much work on general DFO \cite{conn1997,kochenderfer2019,conn2009}, there are two notably similar approaches to the trajectory bundle method. The first is the Constrained Optimization by Linear Interpolation (COBYLA) solver, where the cost and constraint functions are approximated with linear interpolation \cite{powell1994}, and the second is the Gauss-Newton method of \cite{cartis2019}, where samples are used for linear interpolation. The trajectory bundle method builds on these two methods and specializes to trajectory optimization problems where the decision variables at different time steps are only coupled via the dynamics constraints. By exploiting this problem-specific structure and massively parallel simulation, the trajectory bundle method is a simple and robust method for solving trajectory optimization problems to tight constraint and optimality tolerances.

\subsection{Trajectory Optimization through Sequential-Convex \\Programming}
{Trajectory optimization provides a rigorous and powerful mathematical framework for solving general optimal control problems as nonlinear programs. The cost and constraint functions can be arbitrary nonlinear functions that describe the task objective, laws of physics, and physical limits that must be enforced~\cite{malyuta2021}. 
Assuming the cost and constraint functions are smooth and differentiable, this potentially nonlinear, non-convex problem can be solved by linearizing around a current iterate, forming a convex approximation of the original problem, and iterating until convergence. This method, known as Sequential Convex Programming (SCP), has been applied successfully to a wide variety of robotic systems such as rockets \cite{blackmore2016}, orbital transfers \cite{thorne}, fixed-wing aircraft \cite{cory2008}, rotorcraft \cite{sun2022, szmuk2019a}, autonomous vehicles \cite{chen2021, lecleach}, and underwater vehicles \cite{Mamakoukas2016sequential, lee2023}.
}

{The trajectory bundle method solves the trajectory optimization problem through sequential-convex programming in a similar fashion to existing work, but differs in the way in which a convex approximation of the problem is formed. Instead of approximating the nonlinear cost and constraint functions by linearizing around the current iterate, we approximate the cost and constraint functions in a derivative-free manner by linearly interpolating samples within a trust region.}

\subsection{Model-Predictive Path Integral Control}
{
Model-predictive path integral control (MPPI) is a sampling-based method for trajectory optimization or model-predictive control in which derivatives of the dynamics and cost functions are not required. The MPPI algorithm was first derived using a path-integral approach \cite{williams2016}then later re-derived from an information theoretic perspective \cite{williams2018}, a stochastic search perspective \cite{wang2021}, as well as through the use of mirror descent in an online learning context \cite{wagener2019}.} 

{Given the rise in parallel computer architectures (GPUs and multi-threaded CPUs), MPPI has become a popular approach for tackling challenging real-time optimal control problems~\cite{williams2016, williams2018, alvarez2024realtime, liDROPDexterousReorientation2024} because the most expensive part of the algorithm --- evaluating sampled rollouts and their corresponding costs --- can be done entirely in parallel. Additionally, MPPI can be extremely general and naturally amenable to black-box models learned from real-world data~\cite{williams2016}. Despite its empirical success, MPPI performs poorly on open-loop unstable systems due to its single-shooting nature, and is unable to directly reason about constraints.}

{In this paper, we provide a new interpretation of MPPI as a sequential convex programming method and as a special case of the trajectory bundle method, and generalize sample-based optimal control to handle open-loop unstable systems and arbitrary constraints.
} 

\begin{figure}
    \centering
    \includegraphics[width=0.9\linewidth]{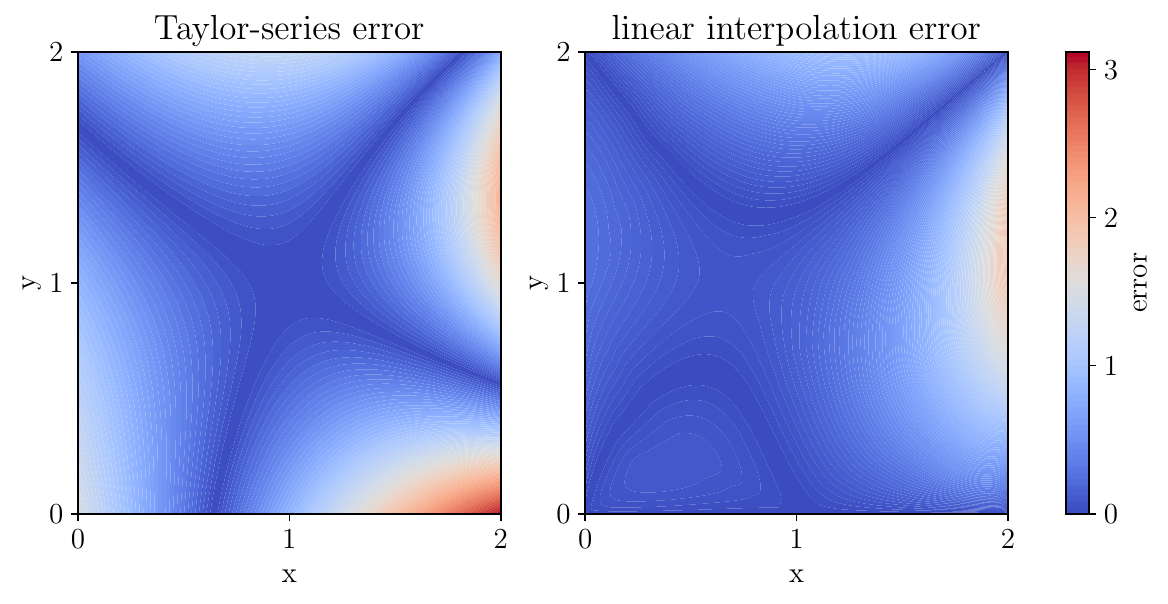}
    \caption{A comparison of the accuracy of a first-order Taylor series taken about $(x,y)=(1,1)$ with linear interpolation of the four corner points on the function $f(x,y) = \sin(x)e^{y}$. While the magnitudes of the errors are comparable between these two approximations, the patterns of these errors are notably different. 
    }
    \label{fig:btb:interp}
\end{figure}
\section{Background}\label{sec:bundles:background}
Like many nonlinear optimization algorithms, the trajectory bundle method handles nonlinear cost and constraint functions by approximating them locally with affine functions. In this section, the standard method of approximating functions with a Taylor series is reviewed, followed by a derivative-free method using linear interpolation over sample points.  Using these linear interpolants, the process for locally approximating a generic constrained optimization problem as convex is then presented. 
\subsection{Affine Function Approximation}\label{btb:sect:interp}
%
%

An arbitrary function $p : \mathbf{R}^a \rightarrow \mathbf{R}^b$ is affine if it can be represented in the following form: $p_\text{aff}(y) = d + Cy$, where $d \in \R{b}$ and $C \in \R{b \times a}$.
The process of locally approximating a nonlinear function with an affine function around a point $\bar{y}$ is often referred to as linearization, with $\bar{y}$ denoted as the linearization point. In this section, the standard method of approximation by first-order Taylor series is presented, followed by a derivative-free method that involves linear interpolation of sampled points.
\subsubsection{Taylor Series}
An affine approximation of this function ${p}(y) \approx \hat{p}(y)$ can be formed in the vicinity of an input value $\bar{y}$ through the use of the first-order Taylor series,
\begin{align}
    p(y) \approx \hat{p}(y) = p(\bar{y}) + \jac{p}{y} (y - \bar{y}),
\end{align}
where both the value and the Jacobian of $p$ are calculated at the point $\bar{y}$. This approximation is exact at $\bar{y}$, and, generally speaking, becomes less accurate farther from $\bar{y}$. 

%
\subsubsection{Linear Interpolation}
Alternatively, nonlinear functions can be approximated in a derivative-free manner by linearly interpolating between sampled function values. This is useful when the derivatives of a function are unavailable, challenging to compute, or unreliable.
%

For an affine function, any linear interpolation between two inputs is equal to the linear interpolation of the outputs. This means for an interpolation parameter $\theta \in [0, 1]$ and two inputs $y_1$ and $y_2$, the following holds:
\begin{align}
    p_\text{aff}(\,\underbrace{\theta y_1 + (1-\theta)y_2}_\text{interpolated inputs}\,) &= \underbrace{\,\theta p_\text{aff}(y_1) + (1 - \theta) p_\text{aff}(y_2)\,}_\text{interpolated outputs}.
\end{align}
This concept can be extended to $m$ points with an interpolation vector $\alpha \in \mathbf{R}^m$ that belongs to a standard simplex:
\begin{align}
    \alpha \in \Delta^{m-1} = \biggl\{  \alpha \in \R{m} \,\,\bigg|\,\, \sum_{i=1}^m \alpha_i =1,\, \alpha \geq 0\biggl\},
\end{align}
where again the convex combination of these $m$ inputs is equal to the same convex combination of the $m$ outputs,
\begin{align}
    p_\text{aff}\bigg({\sum_{i=1}^m\alpha_i y_i}\bigg) &= {\sum_{i=1}^m \alpha_i p_\text{aff}(y_i)}.
\end{align}
This means that we can locally approximate the original nonlinear function $p$ in the neighborhood of $\bar{y}$ by sampling $m$ points from a distribution centered around $\bar{y}$ with $y_i \sim \mathcal{D}(\bar{y})$, and constraining the inputs to this approximation to be a linear combination of the sample points.  For notational convenience, the lists of inputs and outputs are horizontally concatenated as columns of the matrices 
\begin{align}
    W_y &= \begin{bmatrix}
        y_1 & y_2 & \cdots & y_m
    \end{bmatrix} \in \R{n_y \times m}, \\
    W_p &= \begin{bmatrix}
        p(z_1) & p(z_2) & \cdots & p(z_m)
    \end{bmatrix} \in \R{n_r \times m}, 
\end{align}
enabling the affine approximation $\hat{p}$ to be summarized as:
\begin{align}
    y &= W_y \alpha \quad \quad \quad  \text{linear interpolation of inputs,}\\ 
    \hat{p} &= W_p \alpha \quad \quad \quad  \text{linear interpolation of outputs.}\label{btb:blend1}
\end{align}
We are effectively using the approximation $p(W_y\alpha) \approx W_p\alpha$, which is linear in $\alpha$.





\subsection{Approximation for Optimization}
We will now consider a general nonlinear optimization problem and examine how these linearization techniques can be utilized to form a convex approximation of the original problem. This approach is used in sequential convex programming (SCP) methods where nonconvex optimization problems are solved by iteratively approximating the problem as convex in the neighborhood of the local iterate and solving for a step direction \cite{gill2005,nocedal2006,malyuta2021,pantoja1989}. 

To demonstrate how an SCP method works with the two approximation techniques outlined, we examine a generic constrained optimization problem of the following form:
%
\begin{mini}
    {z}{ \| r(z) \|_2^2 }{\label{btb:gen_nl_opt}}{}
    \addConstraint{c(z)}{=0,}{}
\end{mini}
with a decision variable $z \in \R{n_z}$, cost residual function $r : \R{n_z} \rightarrow \R{n_r}$, and constraint function $c : \R{n_z} \rightarrow \R{n_c}$. By approximating both of these functions as affine with a first-order Taylor series around a current iterate $\bar{z}$, we are left with the following convex optimization problem:
\begin{mini}
    {z}{ \| \overbrace{r(\bar{z}) + \frac{\partial r}{\partial z}(z - \bar{z})}^{\hat{r}(z)}\|_2^2 }{\label{btb:lin_prob}}{}
    \addConstraint{\overbrace{c(\bar{z}) + \frac{\partial c}{\partial z}(z - \bar{z})}^{\hat{c}(z)}}{=0.}{}
\end{mini}
While this approximate problem is convex and we are guaranteed to find a globally optimal solution if one exists, we do not have a guarantee of feasibility. There are circumstances in which the linearization of the constraint function results in infeasible approximate problems \cite{nocedal2006}. In order to guarantee that this problem is always feasible, many sequential convex programming methods convert the constraint into a penalty and reformulate \eqref{btb:lin_prob} with an always-feasible variant,
\begin{mini}
    {z, s}{ \| \overbrace{r(\bar{z}) + \frac{\partial r}{\partial z}(z - \bar{z})}^{\hat{r}(z)}\|_2^2 + \mu \|s\|_1}{\label{btb:lin_prob_slack}}{}
    \addConstraint{\overbrace{c(\bar{z}) + \frac{\partial c}{\partial z}(z - \bar{z})}^{\hat{c}(z)} + s}{=0.}{}
\end{mini}
where $\mu \in \R{}_+$ is a positive penalty weight and the $\ell_1$-norm discourages constraint violations.
No matter the structure of the cost and constraint functions, the convex optimization problem in \eqref{btb:lin_prob_slack} is guaranteed to always have a solution. 

Alternatively, a linear interpolant such as that shown in \eqref{btb:blend1} can be used to approximate the cost and constraint functions. To do this, $m$ sample points centered around the current iterate $\bar{z}$ are used to evaluate the cost and constraint functions. These values are then horizontally concatenated into the following matrices:
%
%
\begin{align}
    W_z &= \begin{bmatrix}
        z_1 & z_2 & \cdots & z_m
    \end{bmatrix} \in \R{n_z \times m}, \label{btb:Wz}\\
    W_r &= \begin{bmatrix}
        r(z_1) & r(z_2) & \cdots & r(z_m)
    \end{bmatrix} \in \R{n_r \times m}, \label{btb:Wr}\\
    W_c &= \begin{bmatrix}
        c(z_1) & c(z_2) & \cdots & c(z_m)
    \end{bmatrix} \in \R{n_c \times m}.\label{btb:Wc}
\end{align}
The interpolation vector $\alpha \in \R{m}$ is used to interpolate between these samples and their corresponding cost and constraint values. These approximations are used to recreate the relaxed problem shown in \eqref{btb:lin_prob_slack} as the following:
\begin{mini}
    {\alpha, s}{ \| \overbrace{W_r \alpha}^{\hat{r}(z)}\|_2^2 + \mu \|s\|_1}{\label{btb:lin_prob2}}{}
    \addConstraint{\overbrace{W_c \alpha}^{\hat{c}(z)} + s}{=0}{}
    \addConstraint{\alpha}{\in \Delta^{m-1},}{}
\end{mini}
where the optimal solution is $z^* = W_z \alpha^*$.  It is important to note the similarities between \eqref{btb:lin_prob_slack} and \eqref{btb:lin_prob2}, where the only difference is the method to approximate the cost, residual, and constraint functions. Another key difference between these methods is the implicit trust region present in the simplex constraint on $\alpha$. Since $\alpha \in \Delta ^{m-1}$, the solution is restricted to the convex hull of the sample points. Using a sampling scheme that only samples points within a set trust region, the solution to \eqref{btb:lin_prob2} is guaranteed to stay within the trust region.


\section{The Trajectory Bundle Method}\label{sec:trajectory_bundle_method}
In this section, we outline a canonical trajectory optimization problem specification and use linearly interpolated trajectory bundles to approximate the cost and constraint functions. 

Trajectories will be represented in discrete time as a list of vectors. For a dynamical system with a state $x \in \R{n_x}$ and control $u\in\R{n_u}$, the discrete-time dynamics function $x^{(k+1)} = f(x^{(k)}, u^{(k)})$ maps the state and control at time-step $k$ to the state at $k+1$. A trajectory comprised of $N$ time steps is represented by $(x^{(1:N)}, \, u^{(1:N-1)})$, such that numerical optimization can be used to solve for these values. 

\subsection{Trajectory Optimization}
A canonical trajectory optimization problem considering a trajectory with $N$ time steps is represented as:
\begin{mini}
    {x^{(1:N)}, u^{(1:N-1)}}{  \|r_N(x^{(N)})\|_2^2 + \sum_{k=1}^{N-1}\|r_k(x^{(k)},u^{(k)})\|_2^2}{\label{btb:trajopt}}{}
    \addConstraint{x^{(k+1)}}{=f(x^{(k)}, u^{(k)})}{}
    \addConstraint{c(x^{(k)}, u^{(k)})}{\geq 0,}
\end{mini}
where $c(x^{(k)}, u^{(k)})$ is a generic constraint function. We will assume that all relevant constraints are expressed in this form, including initial and goal constraints, state and control limits, and other general stage-wise constraints.

The problem format in \eqref{btb:trajopt} is often referred to as multiple shooting \cite{hargraves1987,betts2001}, where both the state and control histories are optimized over, and the trajectory only becomes dynamically feasible at convergence. This differs from single shooting, where only the controls are optimized over, and a rollout is performed to recover the states. One important distinction between these two methods is that in single shooting, the discrete-time dynamics must be evaluated sequentially $N-1$ times during the rollout, while in multiple shooting, the $N-1$ dynamics constraints can be evaluated entirely in parallel. This is especially relevant with GPU-based physics simulation, where the speed of a single simulation can be comparable to thousands of simulations run in parallel. 

\subsection{Solving Multiple Shooting with Trajectory Bundles}
The trajectory bundle method is able to reason about the trajectory optimization problem in \eqref{btb:trajopt} without having to differentiate any of the cost, dynamics, or constraint functions. Instead of using derivatives to approximate these functions with their first-order Taylor series, sampled trajectories near the current iterate are used to evaluate these functions for approximation with linear interpolation. This idea is shown in \ref{btb:sect:interp}
Given an initial guess or current iterate $(\bar{x}^{(1:N)}, \bar{u}^{(1:N-1)})$, the costs, constraints, and dynamics functions are computed for each of the $M$ samples surrounding each knot points. To demonstrate this, let us examine a single knot point, $k$, where the current iterate is $(\bar{x}^{(k)}, \bar{u}^{(k)})$. From here, $m$ points are sampled near the iterate, and these samples are horizontally concatenated into the following matrices:
\begin{align}
    W_x^{(k)} &= \begin{bmatrix}
        x^{(k)}_1 & x^{(k)}_2 & \cdots & x^{(k)}_m
    \end{bmatrix} \in \R{n_x \times m}, \label{btb:Wx}\\
    W_u^{(k)} &= \begin{bmatrix}
        u^{(k)}_1 & u^{(k)}_2 & \cdots & u^{(k)}_m
    \end{bmatrix} \in \R{n_u \times m}, \label{btb:Wu}
\end{align}
after which, all of the cost, dynamics, and constraint functions are computed and stored in a similar fashion,
\begin{align}
    W_r^{(k)} &= \begin{bmatrix}
        r_1^{(k)} & r_2^{(k)} & \cdots & r_m^{(k)}
    \end{bmatrix} \in \R{n_r \times m}, \label{btb:Wr2}\\
    W_f^{(k)} &= \begin{bmatrix}
        f_1^{(k)} & f_2^{(k)} & \cdots & f_m^{(k)}
    \end{bmatrix} \in \R{n_x \times m}, \label{btb:Wf2}\\
    W_c^{(k)} &= \begin{bmatrix}
        c_1^{(k)} & c_2^{(k)} & \cdots & c_m^{(k)}
    \end{bmatrix} \in \R{n_c \times m}, \label{btb:Wc2}
\end{align}
where $r_i^{(k)} = r(x^{(k)}_i, u^{(k)}_i)$, $f_i^{(k)} = f(x^{(k)}_i, u^{(k)}_i)$, and $c_i^{(k)} = c(x^{(k)}_i, u^{(k)}_i)$. 
These matrices are computed for time-steps $1\rightarrow N$, with time-step $N$.
Together, these matrices can be used to locally approximate the potentially nonconvex optimization problem in \eqref{btb:trajopt}  as the following convex optimization problem:

\begin{mini}
    {\alpha^{(1:N)}}{ \|\overbrace{W_r^{(N)} \alpha^{(N)}}^{\hat{r}_N(x_N)} \|_2^2 + \sum_{k=1}^{N-1} \| \overbrace{W_r^{(k)} \alpha^{(k)}}^{\hat{r}_k(x_k, u_k)}\|_2^2 }{\label{btb:trajopt_bundled}}{}
    \breakObjective{+ \mu \sum_{k=1}^{N-1}\|s^{(k)}\|_1 + \|w^{(k)}\|_1}
    \addConstraint{\overbrace{W_x^{(k+1)} \alpha^{(k+1)}}^{x^{(k+1)}}}{=\overbrace{W_f^{(k)} \alpha^{(k)}}^{\hat{f}(x^{(k)}, u^{(k)})} + s^{(k)}}{}
    \addConstraint{\overbrace{W_c^{(k)} \alpha^{(k)}}^{\hat{c}(x_k, u_k)} + w_k}{\geq 0}{}
    \addConstraint{\alpha^{(k)}}{\in \Delta^{m-1}\,,}{}
\end{mini}
where $s$ and $w$ are slack variables that, analogous to~\eqref{btb:lin_prob_slack}, ensure the problem remains feasible.

Each time this problem is solved, the new iterates $(x^{(1:N)}, u^{(1:N-1)})$ are used to generate $m$ new samples, and the problem is formed and solved again. This SCP-based algorithm repeats until convergence which, in this particular case, is synonymous with constraint satisfaction.

\subsection{MPPI as a Trajectory Bundle Problem}\label{subsec:mppi}



In MPPI, control policies are sampled and used to generate simulated rollouts with associated costs, and a weighted average based on these costs is used to ``blend'' the sampled policies, with this process at each controller call, the nominal control policy is converging to local optimality.

A single iteration of MPPI starts with $m$ sampled control policies $U^{(1:m)}$ from a distribution centered around a nominal policy $\bar{U}$. Each of these samples is used in a forward dynamics rollout from an initial condition $x_0$ and an associated cost $J_i$ is computed using the rollout from sample $i$. Using the costs from these rollouts, $J_{1:m} \in \R{m}$, weights $\alpha \in \R{m}$ are computed with the softmax function:
\begin{align}
    \alpha_i = \frac{e^{-\frac{J_i}{\lambda}}}{\sum_{i=1}^m e^{-\frac{J_i}{\lambda}}}, \label{softmax_solution}
\end{align}
{where $\lambda$ is a non-negative temperature parameter.}
Using these weights, the resulting updated control policy is a weighted average of the samples, computed as 
\begin{align}
    U = \sum_{i=1}^m \alpha_i U^{(i)}
\end{align}
{In the limit $\lambda \rightarrow 0$, the MPPI update selects the single best control sequence among the samples, which is also known as predictive sampling ~\cite{Howell2022predictive}.}

The MPPI update rule can also be derived as a special case of the trajectory bundle method, where the trajectory is represented with single shooting and there are, therefore, no explicit dynamics constraints.
 Given the $m$ control samples and associated costs, a convex optimization problem solving for the convex combination of trajectories that minimizes the interpolated cost along with a negative entropy regularizer is the following:
\begin{mini}
    {\alpha}{ \overbrace{w_J^T\alpha}^{\hat{J}(u)}  - \overbrace{\lambda \sum_{i=1}^{m} \alpha_i \log \alpha_i}^{\text{entropy regularization}}}{\label{btb:mppi_bundle}}{}
    \addConstraint{\alpha}{\in \Delta^{m-1}.}{}
\end{mini}
where $\lambda$ is the regularization parameter, and the convexity of the negative entropy term makes this a convex optimization problem. The solution to this problem can be computed in closed form as \eqref{softmax_solution}. Just like in MPPI, $\lambda \rightarrow 0$ corresponds to the unregularized bundle problem, where the solution is simply the best sample as no linear combination of the samples can produce a lower cost.
This interpretation of MPPI gives a new perspective through the lens of convex optimization, which can both provide a deeper understanding of the algorithm and provide opportunities for future work. 

\subsection{Sampling Strategies}
{We implement multiple sampling strategies, including Gaussian and uniform distributions, but found that performance was largely independent of the distribution. The experiments in the paper use a simple deterministic coordinate-wise perturbation scheme:
\begin{align}
    z_i &= \bar{z} + e_i \Delta z_i \quad &\forall i \in [1, n] \\  
    z_{i+n} &= \bar{z} - e_i\Delta z_i \quad &\forall i \in [1, n]  
\end{align}
where $\bar{z}$ is the current iterate, $\Delta z$ defines the trust region, and $e_i$ is the unit vector in the $i$-th coordinate direction, and $z_{2n + 1} = \bar{z}$. We leave detailed analysis on sampling strategies for future work.
}
\section{Experiments and Results}\label{sec:experiments}
%
%
%

{In this section, we solve several example motion planning and control problems with the trajectory bundle method. In particular, we show that, in \emph{a unified framework}, the trajectory bundle method solves challenging problems typically only associated with \emph{either} derivative-based SCP (Sec. \ref{subsec:di}, \ref{subsec:quadrotor}, \ref{subsection:race_car_trajopt}) \emph{or} derivative-free sampling-based (Sec. \ref{subsec:race_car_mppi}) methods. Additionally, we show that the trajectory bundle method combines the strengths of both classes of algorithms by solving a class of problem neither method alone can solve (Sec. \ref{subsec:neural_cartpole}). }  

{Unless otherwise noted, the convex optimization problems are solved with Clarabel\footnote{For the examples in~\ref{subsec:race_car}, we temporarily used GUROBI~\cite{achterberg2019s} due to a CVXPY bug in transcribing large problems to Clarabel. We are working with the CVXPY developers to resolve this issue for the camera-ready version.} \cite{goulart2024} through CVXPY \cite{diamond2016cvxpy} and we consider TBM to be converged when the maximum constraint violation reaches below $10^{-4}$.  An open-source implementation of the solver and experiments will be made available upon publication. 

\begin{table}[h]
    \centering
    \begin{tabular}{l|c|c|c}
        \toprule
        \textbf{Features} & \textbf{SCP} & \textbf{MPPI} & \textbf{TBM (ours)} \\
        \midrule
        State constraints & \cmark & \xmark & \cmark \\
        Reasoning over long horizons & \cmark & \xmark & \cmark \\
        Non-differentiable costs and constraints & \xmark & \cmark & \cmark \\
        Black-box models & \xmark & \cmark & \cmark \\
        Unstable systems & \cmark & \xmark & \cmark \\
        \bottomrule
    \end{tabular}
    \caption{Comparison of different trajectory optimization methods.}
    \label{tab:planning_comparison}
\end{table}
\begin{figure}
    \centering
    \includegraphics[width=0.9\linewidth]{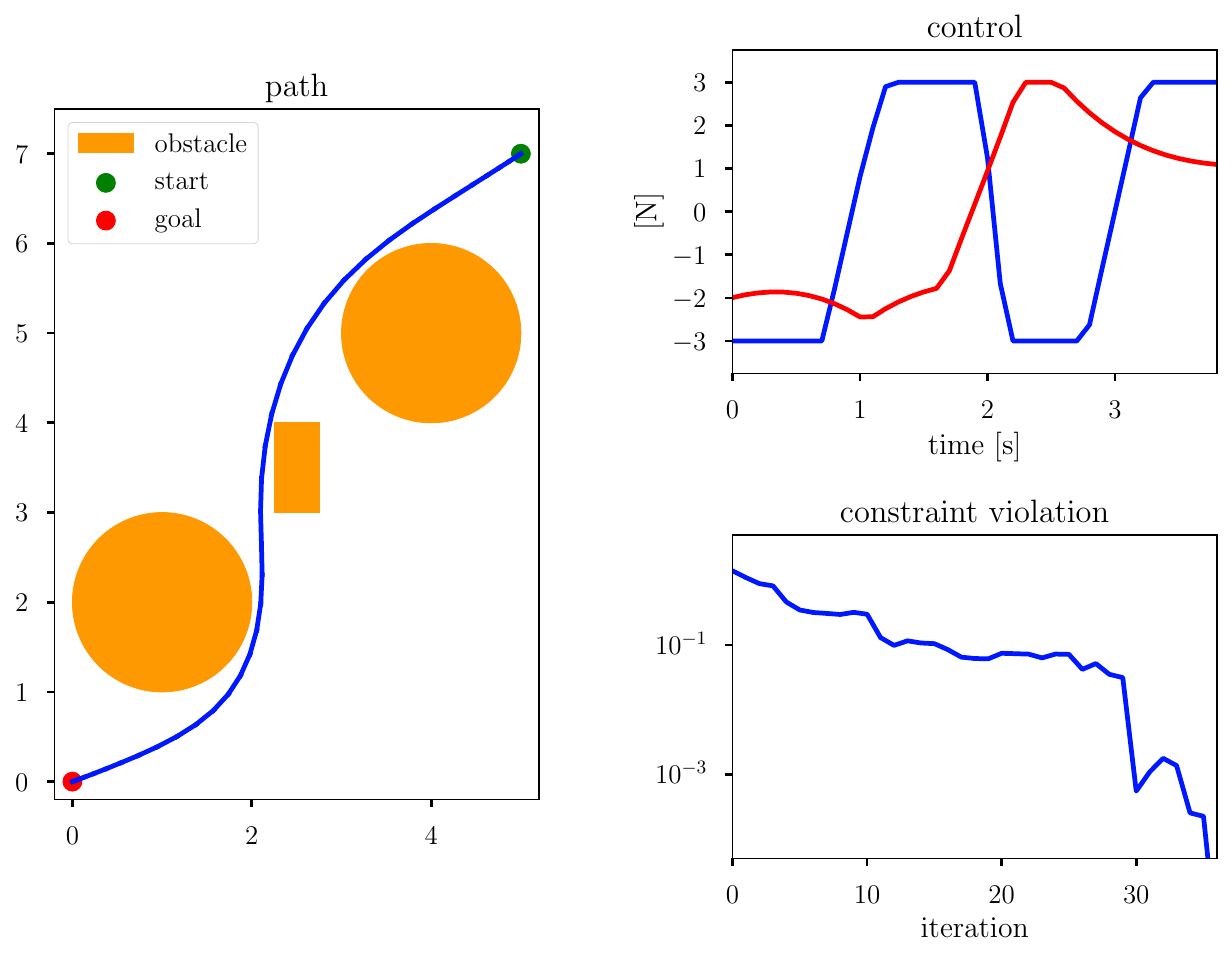}
    \caption{A double integrator with acceleration control is tasked with navigating around three obstacles to a goal position. The trajectory bundle method can directly reason about these nonlinear, non-convex constraints without derivatives, with strong constraint satisfaction and optimality achieved in fewer than 40 iterations.}
    \label{fig:btb:obstacle}
\end{figure}

{\subsection{Double Integrator Collision Avoidance}\label{subsec:di}}
To demonstrate the ability of the trajectory bundle method to handle nonlinear/nonconvex constraints, a collision avoidance example is shown in Fig. \ref{fig:btb:obstacle}. In this scenario, an acceleration-limited double integrator ($u = \ddot{x}$) must find a collision-free path to the goal.
Without derivative information from these nonconvex constraints, the trajectory bundle method converges to a feasible collision-free trajectory in fewer than 40 iterations.

\begin{figure}
    \centering
    \includegraphics[trim={800 350pt 300 200pt}, clip, width=0.49\linewidth]{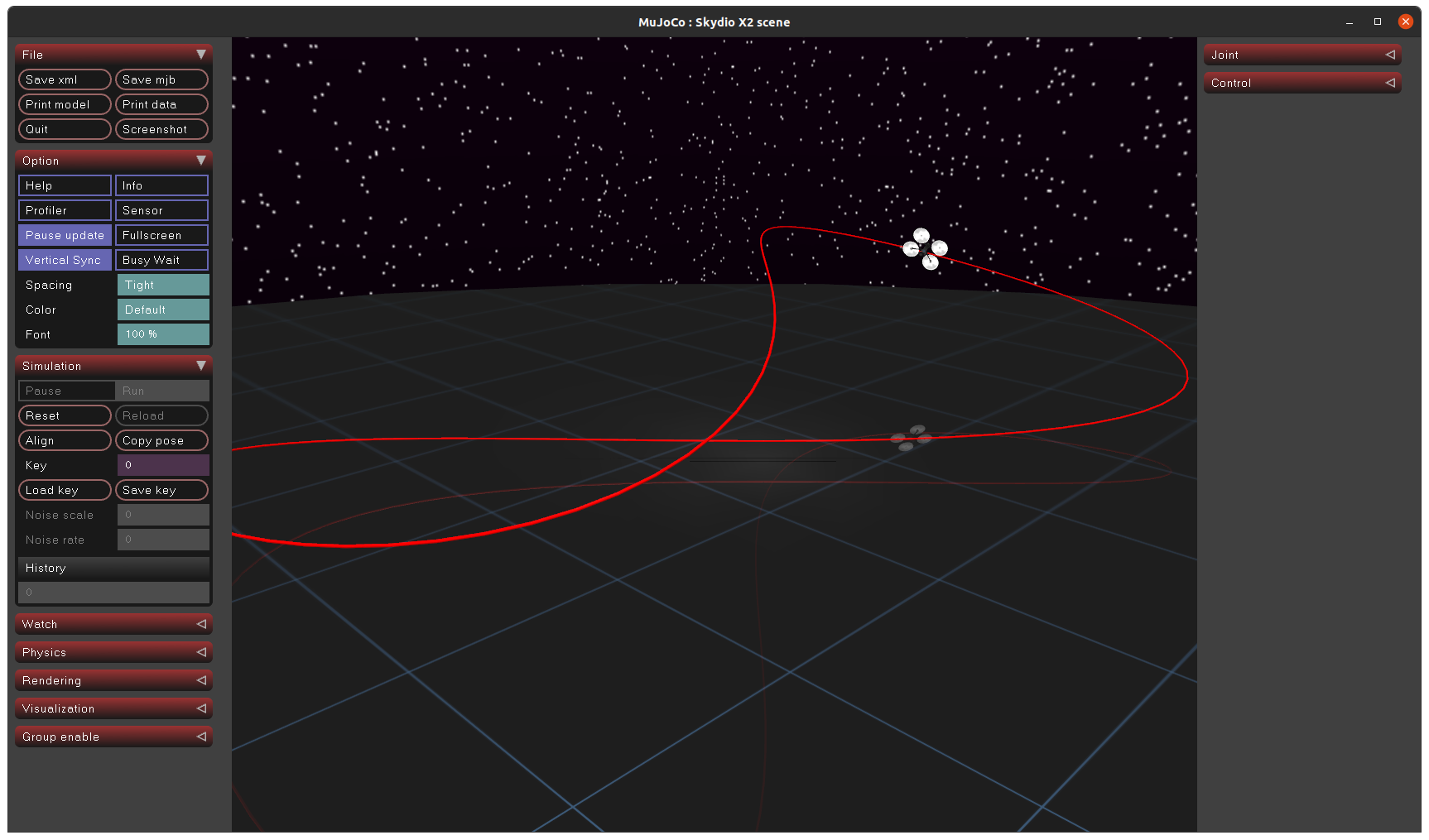}
    \hfill
    \includegraphics[trim={800 350pt 300 200pt}, clip, width=0.49\linewidth]{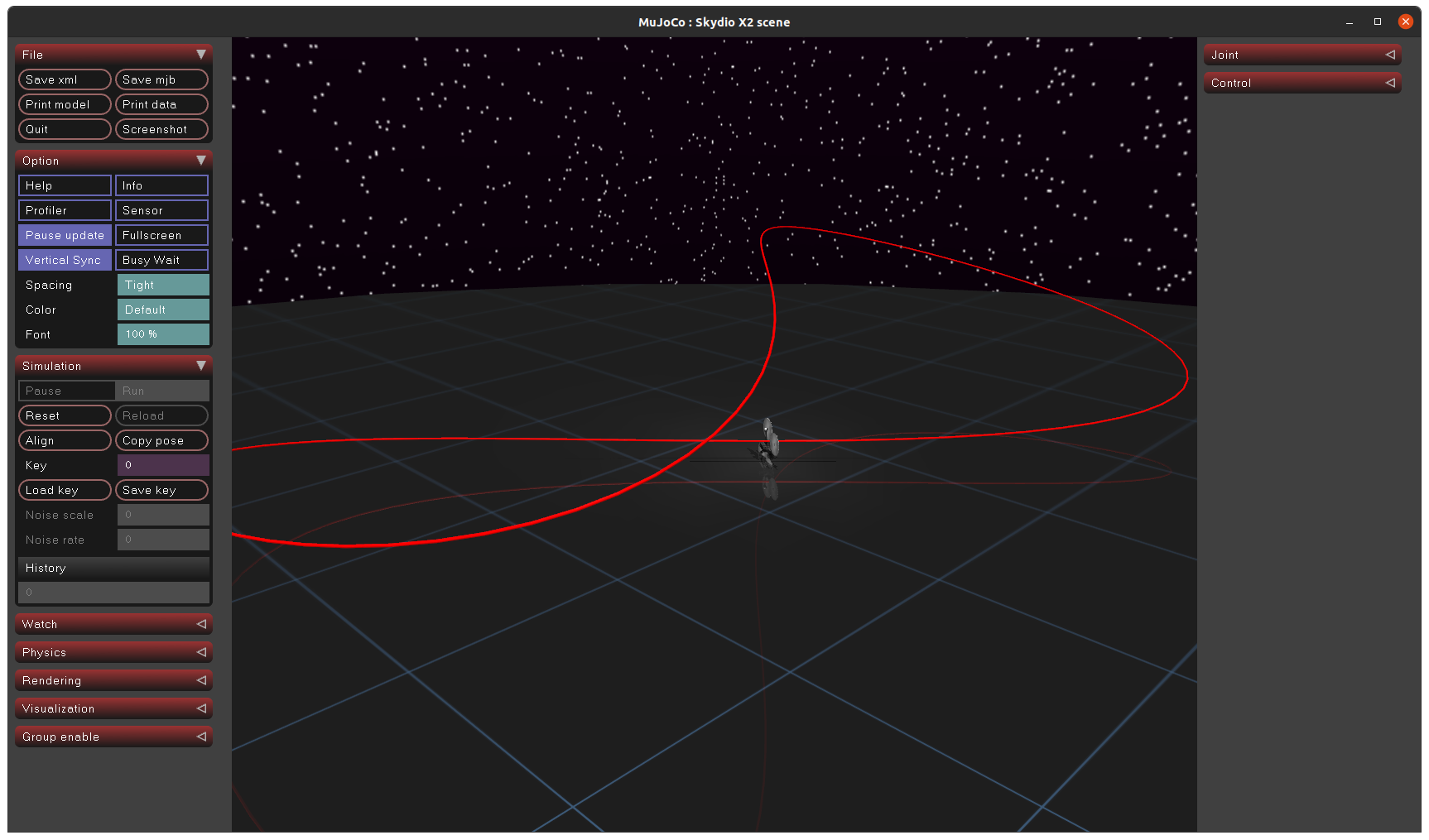}
    \caption{{A quadrotor with rotor-velocity control tracks a figure eight reference over a $5$-second horizon discretized with $100$ time steps. The Trajectory Bundle Method solves for an aggressive and smooth trajectory (left) while MPPI fails (right).}}
    \label{fig:btb:drone_shot}
\end{figure}
{\subsection{Quatrotor Figure Eight Tracking}\label{subsec:quadrotor}}
In Fig. \ref{fig:btb:drone_shot}, a quadrotor with rotor-velocity control is tasked with tracking a skewed figure eight path through space over a five-second horizon. The trajectory is discretized into 100 time-steps, and the resulting optimal trajectory smoothly tracks this aggressive reference while maintaining a smooth control commands. The angular velocity of the quadrotor can reach over 200 degrees per second, where the attitude dynamics are highly nonlinear. In fewer than 60 iterations, the trajectory bundle method can solve this problem to the given constraint tolerance. {In comparison, even after millions of simulation steps, MPPI fails to stabilize this open-loop unstable system over a long horizon, resulting in an unrecoverable crash.
}

{The problems in Sec. \ref{subsec:di} and \ref{subsec:quadrotor} are highly nonlinear/non-convex, highly constrained, and deal with long trajectories. These all present challenges to single-shooting methods like MPPI due to unstable rollouts, a problem noticeably absent from multiple shooting formulations.
\begin{figure}
    \centering
    \includegraphics[width=0.99\linewidth]{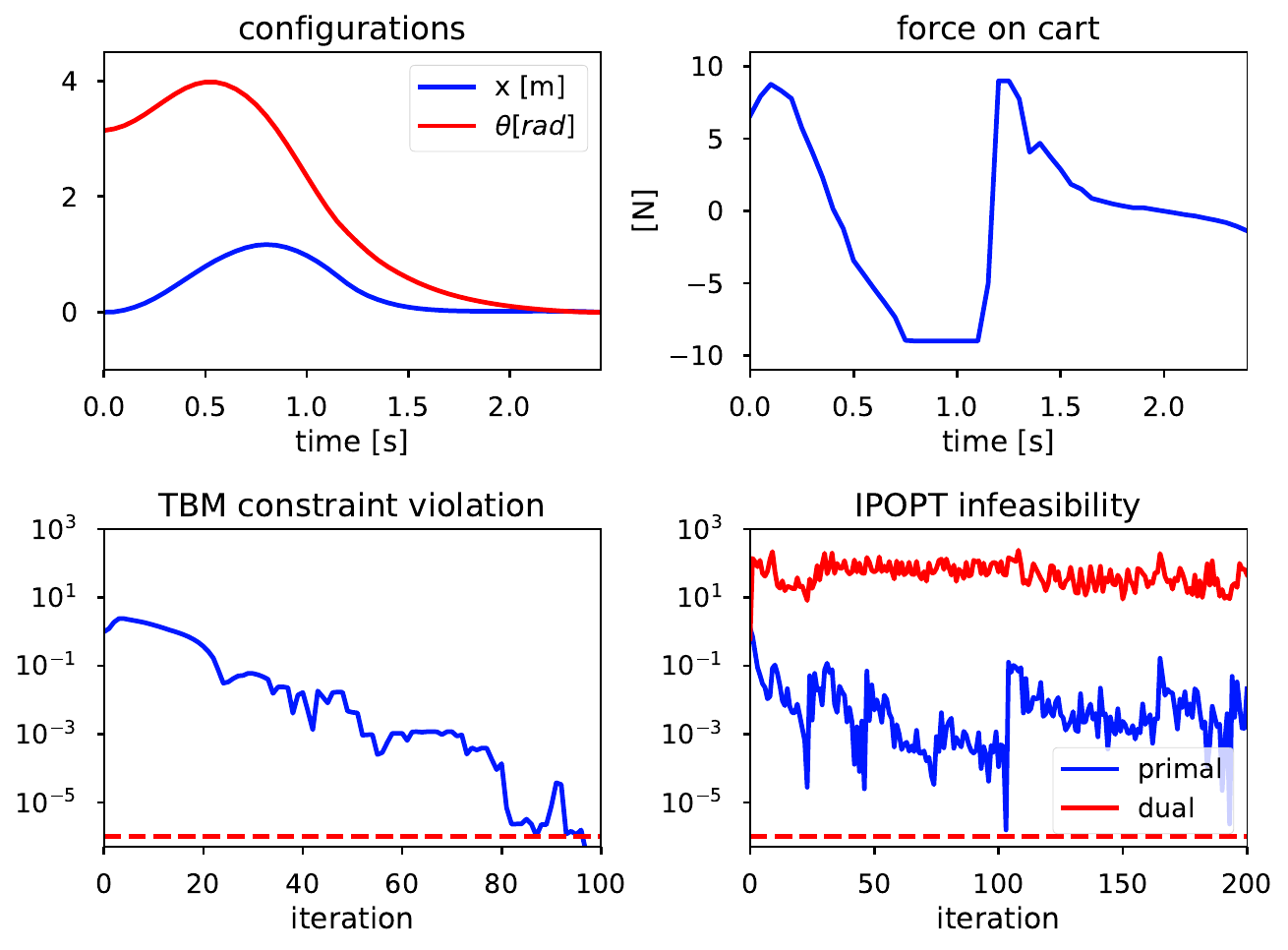}
    \caption{{The cartpole swingup task with a neural-network dynamics model. TBM (ours) finds smooth state (top left) and control (top right) trajectories and converges to tight ($10^{-6}$, red dashed line) constraint satisfaction (bottom left). The baseline IPOPT solver fails to converge (bottom right).}
    }
    \label{fig:btb:cartpole}
\end{figure}
{

\subsection{Cartpole Swingup with a Neural Dynamics Model}\label{subsec:neural_cartpole}
}
{Many robotics simulators are unable to produce smooth and reliable derivatives; this can be a result of nonsmooth impact events, but also if the dynamics are represented with a nonsmooth neural network. In the latter case, it is not uncommon for a learned dynamics model to match the values of the real model well, but not the derivatives. In Fig. \ref{fig:btb:cartpole}, we solve the canonical cartpole swingup trajectory optimization problem with a neural network dynamics model. In this problem, a horizontal cart with an attached pole must swing itself from its stable equilibrium at the bottom to the unstable equilibrium at the top in 2.5 seconds. The problem is discretized into $50$ time steps. Control bounds and goal constraints are also applied.
Using simulated MuJoCo~\cite{todorov2012} dynamics data, we train a 2-layer multilayer perceptron (MLP) with 64 units per layer and ReLU activation functions to predict the next robot state $x^{(k+1)}$ given the current state $x^{(k)}$ and control $u^{(k)}$.
%
Despite the discontinuous MLP dynamics, TBM successfully solves the problem to tight ($10^{-6}$) tolerance in under $100$ iterations. On the other hand, IPOPT~\cite{wachter2006}, a standard nonlinear optimization solver relying on derivative information, fails to converge to an optimal solution even after $10,000$ iterations due to the nonsmoothness of the ReLU network. We only show the first $200$ iterations in Fig. \ref{fig:btb:cartpole} for visual clarity.}

{While sample-based MPPI naturally incorporates potentially nonsmooth learned black-box models, it cannot reason over long horizons or about the general state and goal constraints required in many robotics problems like those described in Sec. \ref{subsec:di}, \ref{subsec:quadrotor}, and \ref{subsec:neural_cartpole}. As a derivative-free trajectory optimization method, TBM is also amendable to these black-box dynamics models and can handle general (potentially black-box) costs and constraints.}

\subsection{Race Car Min-Time Optimization and Local Planning}\label{subsec:race_car}
\begin{figure*}
    \includegraphics[width=\textwidth]{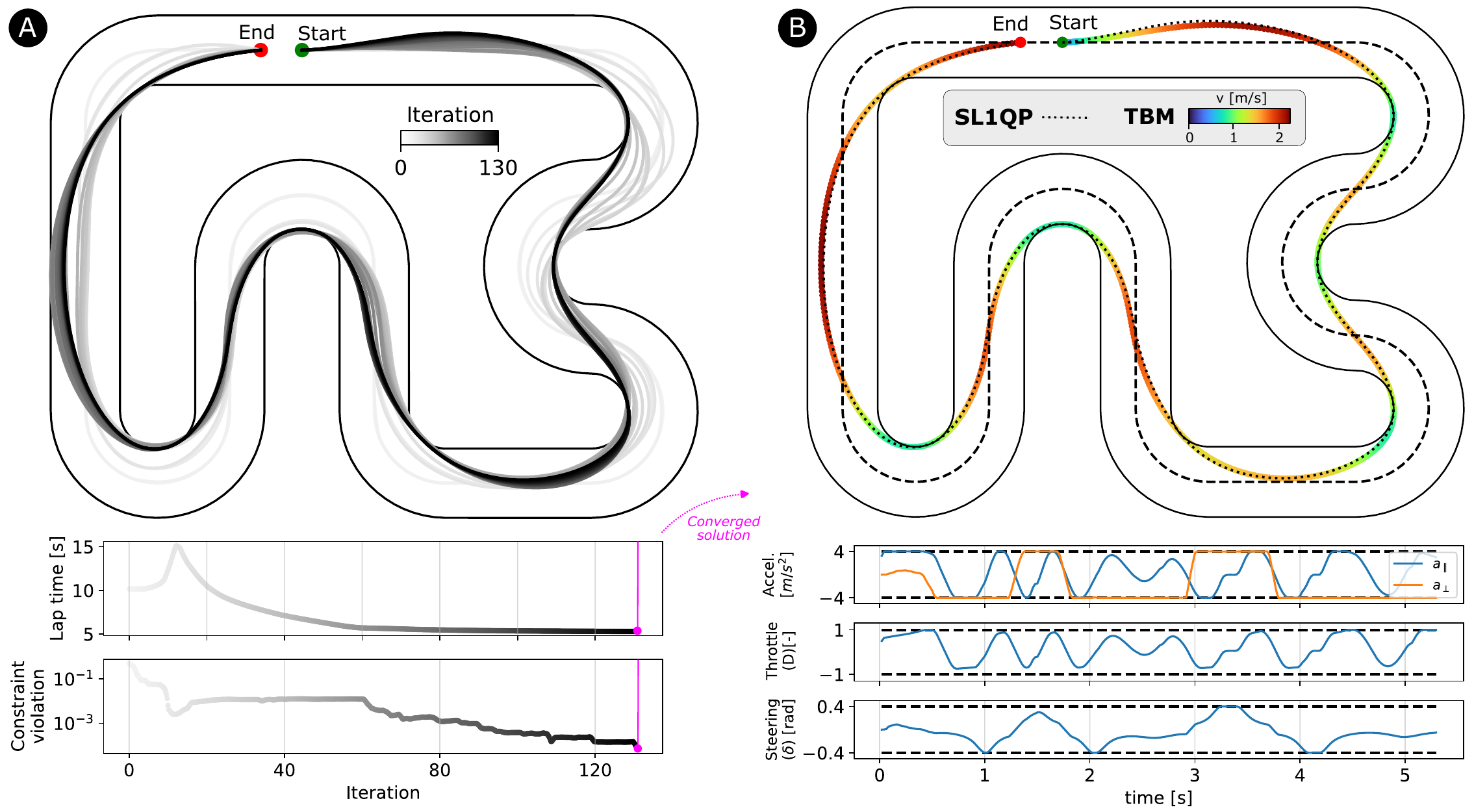}
    \caption{ Minimum-time trajectories using TBM with multiple shooting for a 1:43 autonomous racing example. This demonstrates that TBM (velocity color-coded) and SL1QP (a standard SCP gradient-based method, dotted) produce identical trajectories. \textit{Panel A}: Trajectories progress from centerline initialization (iteration 0) to convergence (iteration 130). The upper section displays the racetrack and trajectories; the lower rows show lap times and constraint violations versus iteration, with grayscale gradients indicating progress. \textit{Panel B}: Converged trajectory at iteration 130. The upper plot shows that the TBM and the SL1QP racelines are the same; the lower section depicts TBM's lateral/longitudinal acceleration, throttle, and steering over time, with dashed lines indicated bounds.}
    \label{fig:car_trajopt}
\end{figure*}

To demonstrate the capability of TBM in generating both SCP and MPPI solutions in a more complex and highly constrained dynamical system, we focus on a state-of-the-art 1:43 scale autonomous car racing example. This case study has been extensively explored in prior research \cite{kloeser2020nmpc,arrizabalaga2023sctomp}, and its time-optimality properties are well-established.

Building on this foundation, we adopt a similar approach to \cite{kloeser2020nmpc} and \cite{arrizabalaga2023sctomp}, utilizing a slip-free bicycle model to describe the vehicle dynamics. The model is defined by the state vector $x = \left[p_x, p_y, \psi, v, D, \delta\right] \in \mathbb{R}^6$ and the input vector $u = \left[\dot{D}, \dot{\delta}\right] \in \mathbb{R}^2$, where $\{p_x, p_y, \psi, v, D, \delta\} \in \mathbb{R}$ represent the car's position, yaw angle, longitudinal velocity, throttle, and steering angle, respectively. For the full equations of motion and model coefficients, we refer the reader to eq.~(3) in~\cite{kloeser2020nmpc}.

To ensure the validity of the system dynamics, we impose constraints on the throttle, steering angle, and their respective time derivatives. Furthermore, to uphold the slip-free assumption---which neglects lateral forces acting on the car---we enforce constraints on the longitudinal and lateral accelerations $\{a_\parallel, a_\perp\} \in \mathbb{R}$. Unlike the state and input constraints, the acceleration constraints are nonlinear functions of the states. All relevant numerical values for these constraints, along with their classification by type, are summarized in Table~\ref{tab:race_car_constraints}.

\begin{table}[t!]
    \renewcommand{\arraystretch}{1.25} 
	\centering
	\caption{Constraints for the race car.} \label{tab:race_car_constraints}
	\begin{tabular}{|c||c|c|c|}
		\hline
		  Type & Constraints\\
        \hline
		  States, $x\in\mathcal{X}$ & $D\in[-1,1]$\,,\;$\delta\in[-0.4,0.4]$\\
		\hline
		Inputs, $u\in\mathcal{U}$  & $\dot{D}\in[-10,10]\,,\;\dot{\delta}\in[-2,2]$\\
        \hline
		  Nonlinear, $g(x)\in\mathcal{G}$  & $\{a_\parallel, a_\perp\}\in[-4,4]$\\
		\hline
	\end{tabular}
\end{table}

Having defined the model, we proceed to demonstrate how TBM is capable of retrieving both SCP and MPPI solutions. To achieve this, we divide our experimental analysis into two parts. First, we focus on a trajectory optimization problem, which is typically solved using SCP. Second, we introduce obstacles to the racetrack to formulate a local planning problem, a scenario commonly addressed with MPPI. Through these two experiments, we aim to show that TBM can match the solutions obtained by both the aforementioned methods in highly constrained and complex systems.

\subsubsection{Trajectory optimization: An SCP case-study}\label{subsection:race_car_trajopt}
In this section, we focus on offline \emph{trajectory optimization}, where the objective is to compute the minimum-time lap. Originally, this problem is an infinite horizon problem, meaning the integration interval over which the optimization is solved depends on the decision variables. Additionally, the solution to this problem is a trajectory that optimally balances the car dynamics and the track constraints. This involves pushing the throttle and steering to the limits of the tires—reaching the edge of the friction circle where the tires begin to slide—while ensuring the car stays within the track boundaries.

Due to these factors, the problem is highly non-convex and nonlinear. The success of local gradient-based methods are reliant on careful initialization of the problem and the use of problem-specific heuristics. These heuristics are necessary to ensure that the validity of the gradients is not compromised by taking excessively large steps. In contrast, TBM is gradient-free and does not require such heuristics, making it easier to implement even in the context of such a challenging problem.

More specifically, the problem is discretized into \( N = 254 \) time steps. The car is initialized in a very simple manner: it is placed at the center of the race track, standing still, and facing forward. The initial state is fully constrained, meaning the car is positioned at the start of the lap with zero velocity. However, for the final state, only the car's position is constrained, while the remaining states (such as velocity and orientation) are left free. This allows the optimization to exploit these degrees of freedom to minimize the lap time. The standard deviations used to generate samples at each iteration are \( \sigma_x = \left[0.06, 0.06, 0.17, 0.5, 0.5, 0.17, 0.01\right] \) for the state variables and \( \sigma_u = \left[1, 0.5\right] \) for the control inputs. 

Additionally, to explicitly minimize the trajectory's time, the vehicle states are augmented by including the time step $\Delta t$ as an additional variable. This provides direct access to the time variable, enabling the reformulation of the cost function in~\ref{btb:trajopt_bundled} as:
\begin{equation}
\min\: T\approx\sum_{k=1}^{N} \|r_k(\Delta t^{(k)})\|_2 + \mu(||s^{(k)}||_1 + ||w^{(k)}||_1)\,,
\end{equation}
where \( \mu = 1 \times 10^7 \) is the penalty coefficient that ensures the satisfaction of the state, input, and nonlinear constraints listed in Table~\ref{tab:race_car_constraints}.

The obtained trajectory, along with the evolution of the solution across different iterations, is illustrated in Fig.~\ref{fig:car_trajopt}. Specifically, in Panel A (left side of the figure), we demonstrate how the trajectories evolve from their initial configuration along the centerline of the track to the time-optimal trajectory over 130 iterations. At the bottom of these plots, the evolution of the lap time and constraint violation is displayed. As expected, both metrics decrease until convergence is achieved at iteration 130.

To further analyze the converged solution, Panel B (right side of Fig.~\ref{fig:car_trajopt}) provides a detailed view of the obtained trajectory. Consistent with racing scenarios, the time-optimal trajectory approaches corners from the inner side and exits toward the outside. The lower part of Panel B displays the longitudinal and lateral accelerations, as well as the throttle and steering inputs. These plots reveal that the car's actuation consistently operates at its physical limits, with at least one acceleration component remaining saturated throughout most of the trajectory, ensuring a motion profile that effectively minimizes lap time.

To validate the solution, we compared the TBM trajectory against a benchmark solution obtained from SL1QP, a widely-used gradient-based solver~\cite{boiroux2019sequential}. The results, illustrated by the dotted line in Panel~B (upper section), demonstrate remarkable agreement between the two methods. The minimal difference in lap time --just 5ms-- can be attributed to variations in numerical precision and solver tolerances. This close correspondence validates that TBM can generate solutions equivalent to traditional gradient-based methods, even when applied to large-scale problems with complex constraints.


\subsubsection{Local Planning: A MPPI case-study}\label{subsec:race_car_mppi}

As discussed in Section~\ref{subsec:mppi}, TBM yields the same solution as MPPI control when an entropy regularization term is added to the cost function. To illustrate this equivalence, we adapt the earlier example of autonomous racing to a common case study where MPPI is frequently applied. Specifically, we introduce $9$ randomly placed obstacles along the race track. This modification renders the previously computed offline trajectory invalid, necessitating the design of a \emph{local planner} capable of completing a lap while avoiding collisions with the obstacles.

In line with standard MPPI practices and in contrast to traditional gradient-based MPC approaches, predictions are generated by sampling over the solely the control sequence (the derivatives of throttle $\dot{D}$ and steering $\dot{\delta}$). To ensure smoothness and adherence to input bounds, we parameterize these inputs using B-Splines and sample over their control points. This approach leverages the property that every point on a B-Spline curve lies within the convex hull of its control points, thereby guaranteeing that the input constraints are satisfied~\cite{arrizabalaga2022spatial}.

Once the inputs are sampled, they are simulated in parallel rollouts, and each trajectory is evaluated using a cost function designed to maximize progress along the race track. Due to the single-shooting nature of MPPI—where only inputs are sampled—we lack the ability to directly impose state constraints. To address this limitation, we reject any samples that violate the constraints outlined in Table~\ref{tab:race_car_constraints}, collide with obstacles, or deviate from the track boundaries. This ensures that only feasible trajectories are considered when computing the weighted average of the sampled inputs.

For this study, we model the system inputs using B-Splines with $100$ control points, generating $m=200$ samples, each comprising $N=20$ time steps. To facilitate a fair comparison between TBM and MPPI, we implement a single-shooting version of TBM by sampling exclusively over the inputs and incorporating the entropy regularization term~\eqref{btb:mppi_bundle} into the cost function. The temperature parameter for the entropy term is set to $\lambda=1 \times 10^{-7}$.

Following the mathematical derivations in Section~\ref{subsec:mppi}, the Trajectory Bundle Method (TBM) and Model Predictive Path Integral (MPPI) control produce identical trajectories, depicted by the magenta trajectory in Fig.~\ref{fig:car_localplanning}. To provide further insight into this solution, the car's motion is illustrated sequentially using magenta boxes. Each frame also displays the sampled trajectories in gray, along with the TBM/MPPI solution highlighted in red. In Fig.~\ref{fig:car_localplanning}, the car successfully completes the lap while avoiding all the obstacles represented as yellow boxes. This example serves as a clear demonstration of the equivalence between TBM and MPPI, showcasing how TBM can achieve the same results as MPPI in practice.

\begin{figure}[t]
    \includegraphics[width=\columnwidth]{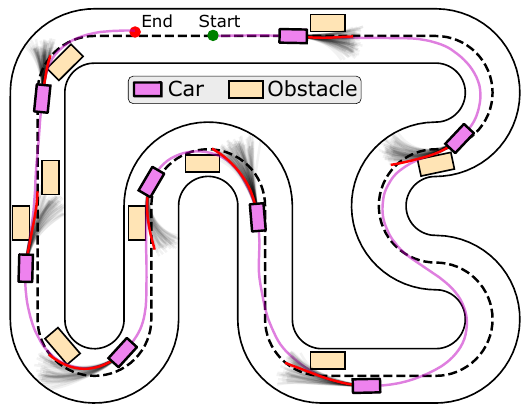}
    \caption{Collision avoidance trajectory using TBM in single shooting fashion for a 1:43 autonomous driving scenario. This example demonstrates that TBM with entropy regularization in~\eqref{btb:mppi_bundle} and MPPI yield identical trajectories. The car's motion is shown with magenta boxes, the obstacles are in yellow, sampled trajectories in gray, and the TBM/MPPI solution in red.}
    \label{fig:car_localplanning}
\end{figure}

\section{Conclusions}\label{sec:conclusions}
In this work, we present the trajectory bundle method, a derivative-free trajectory optimization technique capable of solving nonconvex constrained optimization problems with strong constraint satisfaction. Instead of approximating the nonconvex functions with first-order Taylor series, the trajectory bundle method samples points locally and computes the cost, dynamics, and constraint functions for each of these samples in what we refer to as bundles. These bundles are used to linearly interpolate between these sampled values to approximate the cost, dynamics, and constraint functions. After the computation of these highly parallelizeable function calls, a convex optimization problem is solved where the nonconvex functions are replaced with linear interpolants, and the solution is used to generate new samples for the bundles.  The effectiveness of this method is demonstrated on a variety of robotics platforms.
\section{limitations}
While the trajectory bundle method is a flexible and capable framework for solving trajectory optimization problems, it is not without limitations. Firstly, while you can readily compute constraint violations for these problems, a reliable metric for optimality is still an open question since derivatives are required to compute optimality conditions. Another challenge has been identifying the ideal distribution to draw samples from, while uniform and Gaussian distributions have been sufficient for the examples shown in this paper, there are certainly opportunities to use more expressive and potentially learned distributions for better convergence. The last is a robust convergence rate guarantee, one that will likely require an adaptive penalty $\mu$ and adaptation of the trust region. We leave these challenges open for future work in the area. 

\bibliographystyle{plainnat}
\bibliography{references}

\end{document}